\documentclass[reqno,11pt]{amsproc}
\usepackage{amssymb}
\usepackage{amscd}
\usepackage{amsmath}

\theoremstyle{plain}

\newtheorem{proposition}{Proposition}
\newtheorem{remark}{Remark}

\numberwithin{equation}{section}

\def\supp{\operatorname{supp}}
\def\res{\operatorname{res}}

\def\Nset{\mbox{I\kern-.21em N}}
\def\RE{{\mbox{\rm I\kern-.21em R}}}
\def\ZZ{{\mbox{\sf Z\kern-.45em Z}}}
\def\vv{\kern.344em{\rule[.18ex]{.075em}{1.32ex}}\kern-.344em}

\def\<{\langle} \def\>{\rangle}

\begin{document}

\title{On an inverse problem in photoacoustic.}
\author{M.I. Belishev}%
\address{
    St.Petersburg   Department   of   V.A.Steklov    Institute   of   Mathematics
    of   the   Russian   Academy   of   Sciences, 7, Fontanka, 191023
    St. Petersburg, Russia and Saint Petersburg State University,
    St.Petersburg State University, 7/9 Universitetskaya nab., St.
    Petersburg, 199034 Russia.}\email{ belishev@pdmi.ras.ru  }%
\author{D. Langemann}%
\address{Technische Universit{\"a}t Braunschweig, Inst. Computational
    Mathematics, AG PDE, Universit{\"a}tsplatz 2, 38106 Braunschweig, Germany}
\email{ d.langemann@tu-bs.de  }%
\author{A.S. Mikhaylov} 
\address{St. Petersburg   Department   of   V.A. Steklov    Institute   of   Mathematics
    of   the   Russian   Academy   of   Sciences, 7, Fontanka, 191023
    St. Petersburg, Russia and Saint Petersburg State University,
    St.Petersburg State University, 7/9 Universitetskaya nab., St.
    Petersburg, 199034 Russia.} \email{mikhaylov@pdmi.ras.ru}
\author{V.S. Mikhaylov} 
\address{St.Petersburg   Department   of   V.A.Steklov    Institute   of   Mathematics
    of   the   Russian   Academy   of   Sciences, 7, Fontanka, 191023
    St. Petersburg, Russia and Saint Petersburg State University,
    St.Petersburg State University, 7/9 Universitetskaya nab., St.
    Petersburg, 199034 Russia.} \email{ vsmikhaylov@pdmi.ras.ru }

\date{May, 2018}

\maketitle

\begin{center}
{\bf Abstract.} We consider the problem of reconstruction of the
Cauchy data for the wave equation in $\mathbb{R}^3$ and
$\mathbb{R}^2$ by the measurements of its solution on the boundary
of the unit ball.

\end{center}

\section{Introduction}
Photoacoustic imaging is a technology of visualization for obtaining optical
images in a scattering medium, including  biological tissues \cite{XW}. The basis of
PAI is the effect of laser radiation with ultrashort
pulses on the biological tissue. The absorption of the emitted light pulses by the tissues causes thermoelastic stresses in the light absorption region, a phenomenon known as an optoacoustic or
 photoacoustic effect. Such thermal expansion causes the propagation of ultrasonic
waves in the medium, which in turn can be registered, measured and transformed to produce an image.
The image formation step can be performed by hardware (for example, acoustic or optical focusing)
or by computed tomography, in this case the technology called Photoacoustic tomography.

If we assume that the visualized medium is homogeneous with respect to ultrasound,
we can put the speed of sound equal to $1$ and
write down the standard wave equation for a pressure wave $u(x, t)$, see \cite{KK}
\begin{equation}
    \label{wave_eq}
    \left\{\begin{array}l u_{tt}-\Delta u=0,\quad x\in
        \mathbb{R}^3,\,\,t>0,\\
        u\bigl|_{t=0}=a(x),\,\,u_t\bigl|_{t=0}=b(x)
    \end{array}
    \right.
\end{equation}
The functions $a(x)$, $b(x)$ in the initial conditions describe spatial distributions of the absorption of laser radiation energy and of its temporal change.

Suppose that we have a set of transducers  running over a
surface $S$ to measure the pressure wave and collect data. Suppose that a transducer at
time $t$ located at point $x$ measures  the value of $u(x, t)$. Then we can suppose that
\begin{equation}
\label{observation} F:=u\bigl|_{S\times[0,T]}.
\end{equation}
is known for some $T>0$. The inverse problem (IP) is to find the initial values
$a(x)$, $b(x)$, using the data $F(x,t)$ measured by the  transducers. In
our setting we assume that the surface $S$ is a unit sphere and $\supp a,b$ are located inside $S$.

The first inversion procedures for the spherical case  were
described in \cite{N} in two dimensions  and in \cite{NL} in three
dimensions. The procedure consists of the harmonic decomposition of the
measured data and equating the Fourier coefficients  of the
solutions.

Another procedure was developed in \cite{FPR}, it leads to so
called filtered backprojection type formulas, and it is a most
popular method in photoacoustic at present time.

We propose a new approach, which combines the two listed above and
hopefully leads to a simple and stable numerical scheme.

The article will develop as follows. In the second section we
consider the 3d case, we offer two procedures: in a first one
(exterior) we take the observation $F$ and solve the forward
problem for the wave equation in the exterior of  $B_1\times
[0,T]$, the solution of this problem we denote by $v$. Then we
take the values of $\{v,v_t\}$ at time $t=T$ and solve the wave
equation with this Cauchy data backward in time, the solution of
this problem we denote by $w$. When $T$ satisfies simple geometric
restrictions ($T\geqslant 2$), we get that $a(x)=w(x,0),$
$b(x)=w_t(x,0)$. The second (interior) procedure is based on
solving the problem backward in time.

In the fourth section we study 2d case. The direct application of
the method described in the previous section is impossible due to
the lack of a trailing edge of a wave in 2d. Instead we offer an
iterative procedure.

In the last section we make some remarks on the comparison the
procedure of solving the inverse problem described in the paper
with known ones. We also made the observation that if in
(\ref{wave_eq}) initial displacement or initial speed is zero
(i.e. $a=0$ or $b=0$), then we can recover the remaining part of the
Cauchy data using the observation $u\bigl|_{S\times(0,T)}$ with $T=1$.

The method proposed can be used for solving IP for the wave equation
with known variable speed, that will be a subject of forthcoming publications.
\section{3d case.}

We begin the implementation of the plan proposed in the introduction
with the expansion of the solution of (\ref{wave_eq}) in the
exterior of the unit ball.

\subsection{Solving the IP by "exterior" procedure.}

We are going to construct the solution to the following initial boundary
value problem (IBVP):
\begin{equation}
\label{wave_out} \left\{\begin{array}l
v_{tt}-\Delta v=0,\quad
x\in
\mathbb{R}^3\backslash B_1,\,\,t>0,\\
v\bigl|_{t=0}=0,\,\,v_t\bigl|_{t=0}=0,\\
v\bigl|_{S_1\times [0,T]}=F.
\end{array}
\right.
\end{equation}
We denote $r:=|x|$. Before we pass to the solution to
(\ref{wave_out}), we prove several auxiliary results.

\begin{proposition}
\label{P1} Let $\varphi$ be the solution to
\begin{equation}
\label{wave_vs} \left\{
\begin{array}l
\varphi_{tt}-\Delta \varphi=0,\quad t>|x|,\\
\varphi\bigl|_{|x|=t}=0,
\end{array}
\right.
\end{equation}
such that $\varphi(\lambda x,\lambda t)=\varphi(x,t)$ for $\forall\lambda>0$, then the
function
\begin{equation*}
\psi(x,t)=\int_{0}^{t-r}\omega(\tau)\varphi(x,t-\tau)\,d\tau.
\end{equation*}
is also a solution to (\ref{wave_vs}).
\end{proposition}
\begin{proof}
The condition at $t=r$ are trivially satisfied. Let us evaluate
using that $\varphi(x,r)=0$:
\begin{eqnarray*}
\psi_{x_k}=\int_{0}^{t-r}\omega(\tau)\varphi_{x_k}(x,t-\tau)\,d\tau-\omega(t-r)\varphi(x,r)\frac{x_k}{r},\\
\psi_{x_kx_k}=\int_{0}^{t-r}\omega(\tau)\varphi_{x_kx_k}(x,t-\tau)\,d\tau-\omega(t-r)\varphi_{x_k}(x,r)\frac{x_k}{r},\\
\psi_{t}=\int_{0}^{t-r}\omega(\tau)\varphi_{t}(x,t-\tau)\,d\tau,\\
\psi_{tt}=\int_{0}^{t-r}\omega(\tau)\varphi_{tt}(x,t-\tau)\,d\tau+\omega(t-r)\varphi_t(x,r),
\end{eqnarray*}
from where we get that
\begin{equation}
\label{Psi-eq} \psi_{tt}-\Delta
\psi=\omega(t-r)\left(\varphi_t(x,r)+\sum_{k=1}^3\varphi_{x_k}(x,r)\frac{x_k}{r}\right).
\end{equation}
We differentiate $\varphi(\lambda x,\lambda t)=\varphi(x,t)$ with
respect to $\lambda$ to get
\begin{equation*}
\sum_{k=1}^3 x_k\varphi_{x_k}(\lambda x,\lambda t)+t\varphi_t(\lambda x,\lambda
t)=0.
\end{equation*}
Taking in this equation $\lambda=1$ and $t=r,$ we see that the
right hand side of (\ref{Psi-eq}) is zero, that completes the
proof.
\end{proof}
This statement says that, if we succeed in finding a homogeneous
order zero solution $\varphi$ to (\ref{wave_vs}), we immediately get a
whole family of solutions. We will look for the solutions of
(\ref{wave_vs}) in the form
\begin{equation}
\label{Phi_func}
\varphi(x,t)=\Psi\left(\frac{t}{r}\right)Y^m_n(\theta,\phi),\quad n=1,2,3,\ldots,\quad m=0,\pm 1,\pm 2,\ldots,\pm n,
\end{equation}
where $Y^m_n(\theta,\phi)$ is a spherical function (the restrictions
of the homogeneous harmonic polynomials of degree $n$ to the unit sphere). We
recall that in $\mathbb{R}^3$ the Laplace operator in the spherical coordinates
$(r,\theta,\phi)$ is on the form
\begin{equation*}
\Delta=\frac{1}{r^2}\left(\frac{\partial}{\partial
r}\left(r^2\frac{\partial}{\partial r}\right)+\Delta'\right),
\end{equation*}
where $\Delta'$ is a Beltrami-Laplace operator on the sphere, and the
spherical functions satisfy
\begin{equation*}
-\Delta'Y^m_n=n(n+1)Y^m_n.
\end{equation*}
Below we answer the question how one can choose the function
$\Psi$ in (\ref{Phi_func}).
\begin{proposition}
\label{P2} Let the function $\Psi$ be such that
\begin{equation}
\label{Psi}
\left\{
\begin{array}l
(1-x^2)\Psi''+n(n+1)\Psi=0,\quad x>1,\\
\Psi(1)=0,
\end{array}
\right.
\end{equation}
then the function (\ref{Phi_func}) satisfies the conditions of
Proposition \ref{P1}.
\end{proposition}
\begin{proof}
Let us evaluate:
\begin{eqnarray*}
\varphi_{tt}=\frac{1}{r^2}\Psi''\left(\frac{t}{r}\right)Y^m_n(\theta,\phi),\\
\varphi_r=-\frac{t}{r^2}\Psi'\left(\frac{t}{r}\right)Y^m_n(\theta,\phi),\\
\frac{\partial}{\partial
r}\left(r^2\frac{\partial\varphi}{\partial
r}\right)=\frac{t^2}{r^2}\Psi''\left(\frac{t}{r}\right)Y^m_n(\theta,\phi).
\end{eqnarray*}
Plugging these terms into the wave equation written in the
spherical coordinates, we have:
\begin{eqnarray*}
r^2\Bigl[\varphi_{tt}-\frac{1}{r^2}\left(\frac{\partial}{\partial
r}\left(r^2\varphi_r\right)+\Delta'\varphi\right)\Bigr]=\Psi''\left(\frac{t}{r}\right)Y^m_n(\theta,\phi)-\frac{t^2}{r^2}\Psi''\left(\frac{t}{r}\right)Y^m_n(\theta,\phi)\\
-\Psi\left(\frac{t}{r}\right)\Delta'Y^m_n(\theta,\phi)=
Y^m_n(\theta,\phi)\left(\Psi''\left(\frac{t}{r}\right)-\frac{t^2}{r^2}\Psi''\left(\frac{t}{r}\right)+n(n+1)\Psi\left(\frac{t}{r}\right)\right)\\
=\left[\alpha=\frac{t}{r}\right]=Y^m_n(\theta,\phi)\left((1-\alpha^2)\Psi''\left(\alpha\right)+n(n+1)\Psi\left(\alpha\right)\right).
\end{eqnarray*}
It is clear that, when $\Psi$ satisfies (\ref{Psi}), then the right
hand side of the above equality is zero, and $\varphi$ defined in
(\ref{Phi_func}) satisfies wave equation in (\ref{wave_vs}). The
condition $\Psi(1)=0$ yields the boundary condition
$\varphi(x,|x|)=0$ in (\ref{wave_vs}).
\end{proof}
Let $P_n$ be the Legendre polynomials, i.e. the solutions to
\begin{equation}
\label{Legendre}
\left(\left(1-y^2\right)P_n'(y)\right)'+n(n+1)P_n(y)=0,\quad y>1,
\quad n=1,2,3,\ldots..
\end{equation}
\begin{proposition}
\label{P3} The function
\begin{equation}
\label{Q} Q_{n+1}(x)=\int_1^x P_n(y)\,dy
\end{equation}
satisfies (\ref{Psi}).
\end{proposition}
\begin{proof}
We integrate (\ref{Legendre}) from $1$ to $x$ and by use of the
notation (\ref{Q}) we obtain:
\begin{equation}
\label{Q1} (1-y^2)P_n'(y)\Bigl|_{y=1}^{y=x}+n(n+1)Q_{n+1}(x)=0.
\end{equation}
From (\ref{Q}) we have that $Q_{n+1}''(x)=P_n'(x)$, which together
with (\ref{Q1}) leads to
\begin{equation*}
(1-x^2)Q_{n+1}''(x)+n(n+1)Q_{n+1}(x)=0.
\end{equation*}
The fact that $Q_{n+1}(1)=0$ follows from the definition.
\end{proof}
From Propositions \ref{P1}--\ref{P3} we deduce that the function
\begin{equation*}
\varphi(x,t)=Y^m_n(\theta,\phi)\int_0^{t-r}\omega(\tau)Q_{n+1}\left(\frac{t-\tau}{r}\right)\,d\tau
\end{equation*}
with arbitrary smooth $\omega$ is a solution to (\ref{wave_vs}).

We will  look for the solution to (\ref{wave_out}) in a form
\begin{equation}
\label{V_sol} v(x,t)=\left\{\begin{array}l \sum_{n,m}
Y^m_n(\theta,\phi)\int_0^{t+1-r}\omega^m_n(\tau)Q_{n+1}\left(\frac{t+1-\tau}{r}\right)\,d\tau,\,\,
t\geqslant r-1,\\
0,\quad t<r-1.
\end{array}
\right.
\end{equation}
To satisfy a boundary condition at $r=1$, we wright down the
expansion of inverse data in a form:
\begin{equation}
\label{BC_expan} F(\theta,\phi)=\sum_{n,m}
F^m_n(t)Y^m_n(\theta,\phi), \quad n=1,2,3,\ldots,\quad m=0,\pm 1,\pm 2,\ldots,\pm n,
\end{equation}
and equate:
\begin{equation*}
\sum_{n,m}
Y^m_n(\theta,\phi)\int_0^{t}\omega^m_n(\tau)Q_{n+1}(t+1-\tau)\,d\tau=\sum_{n,m}
F^m_n(t)Y^m_n(\theta,\phi).
\end{equation*}
The above equality leads to the following equations on $\omega^m_n$:
\begin{equation*}
\int_0^{t}\omega^m_n(\tau)Q_{n+1}(t+1-\tau)\,d\tau=F^m_n(t).
\end{equation*}
Differentiating these equations and taking into account (\ref{Q})
and condition $Q_{n+1}(1)=0$, we get relations
\begin{equation*}
\int_0^{t}\omega^m_n(\tau)P_{n}(t+1-\tau)\,d\tau=(F^m_n)'(t).
\end{equation*}
Differentiating one more time and counting that $P_n(1)=1$, we
arrive at Volterra-type equations of a second kind on $\omega_n$:
\begin{equation}
\label{Volterra_eqn}
\omega^m_n(t)+\int_0^{t}\omega^m_n(\tau)P_{n}'(t+1-\tau)\,d\tau=(F^m_n)''(t).
\end{equation}
Finally we arrive at
\begin{proposition}
The solution of (\ref{wave_out}) can be constructed in the form
(\ref{V_sol}), where $\omega^m_n$ are solutions to
the Volterra equation (\ref{Volterra_eqn}), with $F^m_n$ being the
Fourier coefficients in the expansion (\ref{BC_expan}).
\end{proposition}

Since the the kernel in (\ref{Volterra_eqn}) depends on the
difference of the arguments, the solution of (\ref{Volterra_eqn})
can be found in an explicit form. We give a sketch of a proof,
following \cite{S_b5}

By hat we denote the standard Laplace transform:
\begin{equation*}
\widehat f(k)=\int_0^\infty f(t)e^{-kt}\,dt.
\end{equation*}
Taking the inverse of (\ref{Volterra_eqn}) we obtain:
\begin{equation}
\label{Volterra_eqn1}
\omega^m_n(t)=(F^m_n)''(t)-\int_0^{t}H_{n}(t-x)(F^m_n)''(x)\,dx,\,\,
n=1,2,\ldots,\,\, m=0,\pm 1,\ldots,\pm n,
\end{equation}
where $(I+\int_0^tP_n\cdot)^{-1}=I-\int_0^tH_n\cdot$. Taking the
Laplace transform of (\ref{Volterra_eqn}) and
(\ref{Volterra_eqn1}) we get for $n=1,2,\ldots$:
\begin{eqnarray*}
\widehat{\omega^m_n}+\widehat{\omega^m_n}\widehat P_n'=(\widehat{ F^m_n)''},\\
\widehat{\omega^m_n}=\widehat{(F^m_n)''}-\widehat H_n \widehat{(F^m_n)''},
\end{eqnarray*}
from where follows
\begin{equation}
\label{H_n} \widehat H_n=1-\frac{1}{1+\widehat{P_n'(t+1)}}, \quad
n=1,2,\ldots
\end{equation}
Using properties of Laplace transform, we evaluate for
$n=1,2,\ldots:$
\begin{equation}
\label{P_n0} \widehat {P_n'(t+1)}=k\widehat {P_n(t+1)}-P_n(1).
\end{equation}
Expanding $P_n(t+1)$ in the Taylor series at $t=0$, we obtain
\begin{equation*}
P_n(t+1)=P_n(1)+P_n'(1)t+P_n''(1)\frac{t^2}{2}+\ldots+P_n^{(n)}(1)\frac{t^n}{n!},\quad
n=1,2,\ldots,
\end{equation*}
using this expansion and properties of Laplace transform, we
deduce that
\begin{equation}
\label{P_n1} \widehat
{P_n'(t+1)}=\frac{1}{k^n}\left(P_n^{(n)}(1)+\ldots+k^{n-2}P_n''(1)+k^{n-1}P_n'(1)\right).
\end{equation}
So from (\ref{H_n}), (\ref{P_n0}) and (\ref{P_n1}) we deduce that
\begin{equation}
\label{H_n1} \widehat
H_n(k)=1-\frac{k^n}{P_n^{(n)}(1)+\ldots+k^{n-2}P_n''(1)+k^{n-1}P_n'(1)+k^n}.
\end{equation}
Let $k_i,$ $i=1,\ldots,n$ be the  roots of the polynomial in
the denominator in (\ref{H_n1}). Taking $\sigma\in R$ to be the
number greater than the real parts of all $k_i$, we take the
inverse Laplace transform:
\begin{eqnarray}
\label{H_n3} H_n(t)=\frac{1}{2\pi
i}\int_{\sigma-i\infty}^{\sigma+i\infty}\widehat
H_n(k)e^{kt}\,dt\\
=\sum_{i=1}^{n}\res_{k=k_i}\frac{-k^ne^{kt}}{P_n^{(n)}(1)+\ldots+k^{n-2}P_n''(1)+k^{n-1}P_n'(1)+k^n}.\notag
\end{eqnarray}
All aforesaid lead to the
\begin{proposition}
The solution to the Volterra equation (\ref{Volterra_eqn}) can be
evaluated by (\ref{Volterra_eqn1}) with the kernel $H_n$ defined
by (\ref{H_n3}).
\end{proposition}

\subsubsection{Reversing the time. Solution of the IP}

We fix the time $T>0$ and introduce the notations:
$V:=v(\cdot,T),$ $V'=v_t(\cdot,T)$, where $v$ is the solution to (\ref{wave_out}). Note that $\supp V,\supp
V'\subset B_{T+1}\backslash B_1.$ We need to construct the
solution to the following Cauchy problem:
\begin{equation}
\label{wave_invert} \left\{\begin{array}l w_{tt}-\Delta w=0,\quad
x\in
\mathbb{R}^3,\,\,t<T,\\
w\bigl|_{t=T}=V,\,\,w_t\bigl|_{t=T}=V'.
\end{array}
\right.
\end{equation}
As is known, the solution to (\ref{wave_invert}) is given by the Kirchhoff
formula:
\begin{equation}
\label{Kirhgoff}
w(x,t)=-\frac{T-t}{|S_{T-t}|}\int_{S_{T-t}(x)}V'\,dS-\frac{d}{dt}\left(\frac{T-t}{|S_{T-t}|}\int_{S_{T-t}(x)}V\,dS\right),
\end{equation}
where $S_r(x)=\{y\in \mathbb{R}^3|\quad |x-y|=r \}$ is the sphere, $|S_r|=4\pi r^2$ its area.
Since our Cauchy data $V,V'$ are supported in $B_{T+1}\backslash B_1$,
and since we are interested in the values of $w$ and $w_t$ at $t=0$ in a
unit ball, we get the natural restriction on $T$: $T>2.$ Then to
recover the Cauchy data for (\ref{wave_eq}) we just take the values
of $w$ and $w_t$ at $t=0$:
\begin{equation}
\label{data_recov} a(x):=w(x,0),\quad b(x):=w_t(x,0).
\end{equation}
We conclude this section by the algorithm of solving the IP by the
"exterior" procedure:
\begin{itemize}
\item[1)] expand the inverse data in a Fourier series over the
spherical functions (\ref{BC_expan}),

\item[2)] look for the solution $v$ to the IBVP (\ref{wave_out})
in the form (\ref{V_sol}),

\item[3)] find the unknown functions $\omega^m_n$, $n=1,2,\ldots,$, $m=0,\pm1\ldots,\pm n$
in (\ref{V_sol}) as a solutions to the Volterra equation of the
second kind (\ref{Volterra_eqn}), or by the explicit integral formula
(\ref{Volterra_eqn1}) with the integral kernel $H_n$ defined in
(\ref{H_n3}),

\item[4)] fix the time $T\geqslant 2$ and set the functions
$V(x):=v(x,T)$, $V'(x):=v_t(x,T)$,

\item[5)] solve the Cauchy problem (\ref{wave_invert}) for the
function $w$ backward in time by the Kirchhoff formula
(\ref{Kirhgoff}),

\item[6)] Recover the Cauchy data of (\ref{wave_eq}) by taking
$w$, $w_t$ at $t=0$, (\ref{data_recov}).
\end{itemize}

The procedure, offered above is local: suppose that we need to
reconstruct initial data at a single point $x\in B_1$. To use the
Kirchhoff formula it is necessary that the base of the cone
$S_T(x)$ does not intersect the sphere of radius $1$. The latter
leads to restriction $|x|+T\geqslant 1$, and the minimal required
time is exactly $T=1+|x|$. This observation yields the following

\begin{remark}
For the reconstruction of the initial data $a(\cdot)$, $b(\cdot)$ of
(\ref{wave_eq}) at the point $x\in B_1$ we need to know the boundary
observation (\ref{observation}) on $S\times (0,T)$ with $T=1+|x|$.
\end{remark}

\subsection{Solving the IP by "interior" procedure.}

Here we describe another procedure of solving the IP, it works
only if the time fulfill $T\geqslant 2$.
We note, that due to the presence of a trailing edge of a wave in
dimension three, the solution to (\ref{wave_eq}) has a property:
if $T>2$ then $u(x,T)=0$, $|x|\leqslant 1$. Therefore we set up
the following IBVP for a wave equation in the interior of a
cylinder in  reverse time:
\begin{equation}
    \label{wave_eq_rt}
    \left\{\begin{array}l u_{tt}-\Delta u=0,\quad
        \quad |x|<1,\ 0<t<2,\\
        u\bigl|_{S_1\times[0,T]}=F,\quad u\bigl|_{t=2}=u_t\bigl|_{t=2}=0.
    \end{array}
    \right.
\end{equation}
Making a change of variables $t$ to $2-t$ we arrive at
\begin{equation}
    \label{wave_eq_0}
    \left\{\begin{array}l u_{tt}-\Delta u=0,\quad
        \quad |x|<1,\ 0<t<2,\\
        u\bigl|_{S_1\times[0,T]}=\widetilde F,\quad u\bigl|_{t=0}=u_t\bigl|_{t=0}=0,
    \end{array}
    \right.
\end{equation}
where $\widetilde F(\cdot,t)=F(\cdot,2-t)$,
and the unknown functions $a(\cdot)$ and $b(\cdot)$ are consequently equal
 $a(x)=u(x,2)$, $b(x)=-u_t(x,2)$. Thus to solve the IP is
to solve the wave equation $(\ref{wave_eq_0})$. It can be done in two
different ways.

\subsubsection{Solution of (\ref{wave_eq_0}) by residues.} We seek for the
solution in a form
\begin{equation}
\label{u_repr} u(x,t)=\sum_{n,m}
Y^m_n(\theta,\varphi)\varphi^m_n(r,t),\quad n=1,2,\ldots,\quad m=0,\pm 1,\ldots\pm n.
\end{equation}
Substituting this expansion into (\ref{wave_eq_0}), and using
properties of spherical functions $Y^m_n$ we obtain that
$\varphi^m_n$ satisfy:
\begin{equation}
    \label{wave_eq_varphi}
    \left\{\begin{array}l (\varphi^m_n)_{tt}-((\varphi^m_n)_{rr}+\frac 2r (\varphi^m_n)_r)+\frac{n(n+1)}{r^2}{\varphi^m_n}=0,\quad
        \quad r<1,\ 0<t<2,\\
        {\varphi^m_n}\bigl|_{r=1}=F^m_n(t),\quad {\varphi^m_n}\bigl|_{t=0}={\varphi^m_n}_t\bigl|_{t=0}=0.
    \end{array}
    \right.
\end{equation}

Taking the Laplace transform of (\ref{wave_eq_varphi}) with
respect to $t$:
$\chi(r,s)=\int_0^{+\infty}e^{-st}\varphi(r,t)\,dt,\,
f^m_n(s)=\int_0^{+\infty}e^{-st}F^m_n(t)\,dt,$ we get
\begin{equation}
    \label{eq_chi}
    \left\{\begin{array}l s^2{\chi^m_n}-((\chi^m_n)_{rr}+\frac 2r (\chi^m_n)_r)+\frac{n(n+1)}{r^2}{\chi^m_n}=0,\\
        {\chi_n}\bigl|_{r=1}=f^m_n(s).
    \end{array}
    \right.
\end{equation}
On introducing the notation  $L(r)=\sqrt r\chi^m_n(r)$ (we omit
indexes $n$ and $m$ here), we get
$$
(\chi^m_n)_r(r)=r^{-1/2}L_r-\frac12r^{-3/2}L,\quad
(\chi^m_n)_{rr}(r)=r^{-1/2}L_{rr}-r^{-3/2}L_r+  \frac34r^{-5/2}L,
$$
and equation (\ref{eq_chi}) is equivalent to the Bessel equation:
\begin{equation*}
{L}_{rr}+\frac 1r {L}_r-
\frac1{r^2}\left((n+1/2)^2+r^2s^2\right)L=0.
\end{equation*}
Therefore $L=c J_{n+1/2}(-irs)$. Using (\ref{eq_chi}) and the
definition of $L$, we arrive at
\begin{equation}
    {\chi^m_n}(r,s)=\frac{J_{n+1/2}(-irs)}{\sqrt r}  \frac{f^m_n(s)}{J_{n+1/2}(-is)}.
\end{equation}
Taking the inverse Laplace transform yields
$$
\varphi^m_n(r,t)=\frac{1}{\sqrt r 2\pi
i}\int_{\sigma-i\infty}^{\sigma+i\infty} \frac
{J_{n+1/2}(-irs)}{J_{n+1/2}(-is)}f^m_n(s)e^{st}\,ds.
$$
The latter expression can be rewritten in the following form
\begin{equation}
\label{ph_n} \varphi^m_n(r,t)=\psi^m_n(r,t) -2  \sum_{p=1}^\infty
\frac{A_p \sin(k_p^nt+\omega_p)}{J_{n-1/2}(k_p^n)} \frac
{J_{n+1/2}(k_pr)}{\sqrt r},
\end{equation}
where $k_p^n$ are zeros of $J_{n+1/2}(\cdot)$,
$f^m_n(ik_p)=A_pe^{i\omega_p}$ and $\psi^m_n(r,t)$ is the sum of
the residues at the poles of $f^m_n$. Note that this formula is valid  if the poles of $f^m_n$ do not coincide with the zeros of the denominator.

\subsubsection{Solution of (\ref{wave_eq_0}) by Volterra equations.}
Instead of computing residues, the functions $\varphi^m_n$,
   can be calculated in the other way. We assume
that they have a form (we omit indexes $n$ and $m$)
\begin{equation}
    \label{form_varphi}
    \varphi(r,t)= \left\{\begin{array}l \frac{1}{r}\int\limits_{t-1-r}\limits^{t-1+r}\omega(\tau)P_n(\frac{\tau+1-t}{r})\,d\tau,\quad t>1-r,\\
        0,\quad t<1-r.
    \end{array}
    \right.
\end{equation}
Changing the variables $\alpha=\frac{\tau+1-t}{r}$ in the
integral in (\ref{form_varphi}) , we receive that
$\varphi(r,t)=\int\limits_{-1}\limits^{1}\omega(r\alpha+t-1)P_n(\alpha)\,d\alpha$. Then we can evaluate:
\begin{eqnarray*}
    \varphi_{tt}-\varphi_{rr}-\frac{2}{r}\varphi_r+\frac{n(n+1)}{r^2}\varphi=\\
    =\int\limits_{-1}\limits^{1}(\omega''(r\alpha+t-1)-\alpha^2\omega''(r\alpha+t-1)-\frac{2\alpha}{r}\omega'(r\alpha+t-1)+\\
    +\frac{n(n+1)}{r^2}\omega(r\alpha+t-1) )P_n(\alpha)\,d\alpha =\\
    =\int\limits_{-1}^1 \Big( \frac{1-\alpha^2}{r^2} \frac{\partial^2}{\partial\alpha^2}\omega(r\alpha+t-1) -
    \frac{2\alpha}{r^2} \frac{\partial}{\partial\alpha}\omega(r\alpha+t-1)+\\
    +\frac{n(n+1)}{r^2}\omega(r\alpha+t-1)
    \Big)P_n(\alpha)\,d\alpha=\\
    =\frac{1}{r^2}\int\limits_{-1}\limits^{1} \omega(r\alpha+t-1)\Big[ \Big((1-\alpha^2)P_n'(\alpha) \Big)'+n(n+1)P_n(\alpha) \Big]\,d\alpha =0.
\end{eqnarray*}
The latter equality is valid if and only if $P_n$ are the Legendre
polynomials. Therefore (\ref{wave_eq_varphi}) is fulfilled if and only if
 $\varphi^m_n(t,1)=F^m_n(t)$, $n=1,2,\ldots$, $m=0,\pm1\ldots,\pm n$, which we can rewrite
as:
$$
\int\limits_{t-2}^{t}
\omega^m_n(\tau)P_n(\tau+1-t)\,d\tau=F^m_n(t).
$$
On differentiating the above equation with respect to $t$ and
using that $P_n(1)=1$, $P_n(-1)=(-1)^n$, we receive the following
relation:
\begin{equation}
\label{Volt_eq_steps}
\omega^m_n(\tau)-(-1)^n\omega^m_n(\tau-2)-\int\limits_{t-2}^{t}\omega^m_n(\tau)P_n'(\tau+1-t)\,d\tau
=(F^m_n)'(t).
\end{equation}
The functions $\omega^m_n(\cdot)$ can be determined step by step, by
solving integral Volterra equations of second kind on intervals
$(2(k-1),2k)$, $k=1,2,\ldots$. On the first two intervals these
equations have form:
\begin{eqnarray*}
\omega^m_n(\tau)-\int\limits_{t-2}^{t}\omega^m_n(\tau)P_n'(\tau+1-t)\,d\tau=(F^m_n)'(t),\quad t\in (0,2)\\
\omega^m_n(\tau)-\int\limits_{2}^{t}\omega^m_n(\tau)P_n'(\tau+1-t)\,d\tau
=(F^m_n)'(t)+(-1)^n\omega^m_n(\tau-2)\\
+\int\limits_{t-2}^{2}\omega^m_n(\tau)P_n'(\tau+1-t)\,d\tau,\quad
t\in (2,4).
\end{eqnarray*}

We conclude this section with the algorithm of solving the IP by
"interior" procedure:
\begin{itemize}
\item[1)] expand the inverse data in Fourier series over the
spherical functions (\ref{BC_expan}),

\item[2)] look for the solution $u$ to the IBVP (\ref{wave_eq_0})
in a form (\ref{u_repr}),

\item[3)] find the unknown functions $\varphi^m_n$, $n=1,2,\ldots$, $m=1,\pm 1,\ldots\pm n$
by exact formula(\ref{ph_n}), or  in a form (\ref{form_varphi}),
where $\omega$ is as a solutions to the integral equation
(\ref{Volt_eq_steps}),

\item[4)] Recover the Cauchy data of (\ref{wave_eq}) by
$a(x)=u(x,2)$, $b(x)=-u_t(x,2)$.
\end{itemize}

\section{2d case.}

In this section, we consider IP for IBVP (\ref{wave_eq}) with
observation (\ref{observation}) assuming that the dimension $n=2$.
A key feature of the two-dimensional case is the absence of a trailing
edge in the spreading wave. In the three-dimensional situation the
solution to (\ref{wave_eq}) with Cauchy data supported in a unit
ball, has a property that $u(x,T)=0$ in $B_1$, when $T$ is big
enough. It is not the case in 2d. So we cannot directly apply the
method from the previous section: we do not know the value of a
solution $u(x,T)$ for $x\in B_1$! We  show that
$u(x,T)$ for $x\in B_1$ is getting smaller as $T$ is getting
larger, and using this fact we can formulate the approximating
procedure.

\subsection{Solving IBVP for 2d wave equation.}

\subsubsection{Solution of IBVP for a wave equation in the exterior of a ball.}

Consider the following IBVP:
\begin{equation}
    \label{wave_out_2} \left\{\begin{array}l
        v_{tt}-\Delta v=0,\quad
        x\in
        \mathbb{R}^2\backslash B_1,\,\,t>0,\\
        v\bigl|_{t=0}=0,\,\,v_t\bigl|_{t=0}=0,\\
        v\bigl|_{S_1\times [0,T]}=F.
    \end{array}
    \right.
\end{equation}
First we note that Proposition \ref{P1} is valid in any dimensions
(here we use it for $n=2$). Therefore, if $\varphi(x,t)$ is a
solution to (\ref{wave_vs}) such that $\varphi(\lambda x,\lambda
t)=\varphi(x,t)$ then for any smooth function $\omega$, the
function
\begin{equation*}
v(x,t)=\int_{0}^{t-r}\omega(\tau)\varphi(x,t-\tau)\,d\tau.
\end{equation*}
is also solution to (\ref{wave_vs}).  As in 3d case, we  seek for
$\varphi(x,t)$ in a form
$$
\varphi(x,t)=\Psi\left(\frac{t}{r}\right)Y^m_n(\phi),
$$
where $r\geqslant 0$ and $0\leqslant \phi\leqslant 2\pi$ are polar
radius and angle. Functions $Y^m_n(\phi)$ are equal to
$Y^1_n(\phi)=\sin(n\phi)$, $Y^2_n(\phi)=\cos(n\phi)$, $n=0,1,2,\ldots$. Plugging such an
ansatz in (\ref{wave_vs}), we see that $\Psi$ should satisfy the
following problem
\begin{equation}
\label{Psi_2}
\left\{
\begin{array}l
(1-x^2)\Psi''(x)-x\Psi'(x)+n^2\Psi(x)=0,\quad x>1,\\
\Psi(1)=0.
\end{array}
\right.
\end{equation}
This equation is an analog of (\ref{Psi}) in dimension $n=2$. On
introducing the function $U$ by the rule:
$$
\Psi(x)=\sqrt{x^2-1}\cdot U(x),
$$
we see that $U$ satisfies
$$
(1-x^2)U''(x)-3xU'(x)+(n-1)(n-1+2)U(x)=0,
$$
which is a differential equations satisfied by Chebyshev
polynomials of the second kind ($U_n$). Therefore $U=U_{n-1}$,
$n=1,2,\dots$ and $\Psi$ has a form
\begin{equation}
\label{form_Psi}
\Psi_n(x)=  \left\{\begin{array}l \sqrt{x^2-1}U_{n-1}(x),\quad n=1,2,\dots\\
  \log(x+\sqrt{x^2-1}) \quad n=0.
\end{array}
\right.
\end{equation}
We will  look for the solution to (\ref{wave_out_2}) in a
form
\begin{equation}
\label{V_sol_2} v(x,t)=\left\{\begin{array}l \sum_{n,m}
Y^m_n(\phi)\int_0^{t+1-r}\omega^m_n(\tau)\Psi_n\left(\frac{t+1-\tau}{r}\right)\,d\tau,\,\,
t\geqslant r-1,\\
0,\quad t<r-1.
\end{array}
\right.
\end{equation}
To satisfy the boundary conditions at $r=1$, we expand $F$
\begin{equation}
\label{BC_expan_2}
F(t,\phi)=\sum_{n,m}^\infty
F^m_n(t)Y^m_n(\phi)
\end{equation}
in spherical functions and equate:
\begin{equation*}
\sum_{n,m}
Y^m_n(\phi)\int_0^{t}\omega^m_n(\tau)\Psi_n(t+1-\tau)\,d\tau=\sum_{n,m}
F^m_n(t)Y^m_n(\phi),
\end{equation*}
which implies the following equations on $\omega_n$:
\begin{equation}\label{Volterra_eqn_2d}
\int_0^{t}\omega^m_n(\tau)\Psi_n(t+1-\tau)\,d\tau=F^m_n(t),\quad
n=0,1,2,\ldots\quad m=1,2.
\end{equation}
On the other hand, differentiating (\ref{Volterra_eqn_2d}) with
respect to $t$ and using the facts that Chebyshev polynomials of
first kind $T_n$ and second kind $U_n$ satisfy the Pell equation
$$
T_n^2(\tau)-(1-\tau^2)U_{n-1}(\tau)=1,\quad n=1,2,\ldots
$$
and
$$
T_n'(\tau)=nU_{n-1}(\tau),\quad n=1,2,\ldots,
$$
reduces (\ref{Volterra_eqn_2d}) to the singular integral equation
\begin{equation}\label{Volterra_eqn1_2d}
\int_0^{t}\omega^m_n(\tau)\frac{nT_n(t+1-\tau)}{\sqrt{(t+1-\tau)^2-1}}\,d\tau=(F^m_n)'(t),\quad
n=1,2,\ldots,\quad m=1,2.
\end{equation}
The unique solvability of (\ref{Volterra_eqn_2d}) is proved in
\cite{S_b4}. Finally we come to the
\begin{proposition}
The solution of (\ref{wave_out_2}) can be constructed in a form
(\ref{V_sol_2}), where $\omega_n$ are found as  solutions to the
singular integral equation (\ref{Volterra_eqn1_2d}), with $f_n$
being Fourier coefficient in the expansion (\ref{BC_expan_2}).
\end{proposition}
As in 3d case we have that for $\quad n=1,2,\ldots$
\begin{equation}
\label{Volterra_eqn1_2}
\omega^m_n(t)=\int_0^{t}H_{n}(t-x)(F^m_n)'(x)\,dx,
\end{equation}
where the kernel $H_n(t)$, $t>0$ is expressed as
\begin{eqnarray}
 H_n(t)=\frac{1}{2\pi
    i}\int_{\sigma-i\infty}^{\sigma+i\infty}\widehat
H_n(k)e^{kt}\,dt
\end{eqnarray}
 and  $\widehat H_n(k)$  connected with kernel of (\ref{Volterra_eqn1_2d}) via
\begin{equation}\label{H_n2}
 \widehat H_n = \frac{1}{\widehat G_n},\quad  G_n(\tau)=\frac{nT_n(t+1-\tau)}{\sqrt{(t+1-\tau)^2-1}}.
\end{equation}
Using the fact that the Laplace transform obeys
$\mathcal{L}(xf(x))=-\mathcal{L}(f)'$, we conclude that the function
$H_n$ can be calculated via functions $T_n(t)$ where instead of
$t^p$ we plug in  $(-1)^pK^{(p)}$, where
$K_0=\mathcal{L}\frac{1}{\sqrt{\tau^2-1}}$ is a modified Bessel
function of a second kind (Macdonald's function).
\begin{proposition}
Functions $\omega^m_n$, $n=1,2,\ldots$ defined by
(\ref{Volterra_eqn1_2}) with the kernel $H_n$ given by
(\ref{H_n2}) are solutions to singular integral equations
(\ref{Volterra_eqn1_2d}).
\end{proposition}

\subsubsection{Solution of IBVP for a wave equation in the interior of
a ball (by residues)}

Note that in this case, due to the lack of a trailing edge in
dimension two, $u(x,T)\not=0$ in $B_1$ for $T>2$! We will look for
a solution to (\ref{wave_eq_rt}) (with $d=2$) in a form
$u(x,t)=Y^m_n(\varphi)\varphi^m_n(r,t)$. Substituting this expression
into the equation, and using properties of the spherical functions
$Y^m_n$, we obtain that $\varphi_n$, $n=1,2,\ldots$ satisfy the
following IBVP:
\begin{equation}
\label{wave_eq_varphi_2}
\left\{\begin{array}l (\varphi^m_n)_{tt}-((\varphi^m_n)_{rr}+\frac 1r (\varphi^m_n)_r)+\frac{n^2}{r^2}{\varphi^m_n}=0,\quad
\quad r<1,\ 0<t<2,\\
{\varphi^m_n}\bigl|_{r=1}=F^m_n(t),\quad {\varphi^m_n}\bigl|_{t=0}=(\varphi^m_n)_t\bigl|_{t=0}=0.
\end{array}
\right.
\end{equation}
Taking the Laplace transform with respect to $t$:
$\chi(r,s)=\int_0^{+\infty}e^{-st}\varphi(r,t)\,dt$, we obtain
that $\chi^m_n$, $n=1,2,\ldots$ satisfy:
\begin{equation}
\label{eq_chi_2}
\left\{\begin{array}l s^2{\chi^m_n}-((\chi^m_n)_{rr}+\frac 1r (\chi^m_n)_r)+\frac{n^2}{r^2}{\chi^m_n}=0,\\
{\chi^m_n}\bigl|_{r=1}=f^m_n(s).
\end{array}
\right.
\end{equation}
Therefore $\chi_n=c J_{n}(-irs)$, $n=1,2,\ldots$. Using boundary
condition in (\ref{eq_chi_2}),  we arrive at
\begin{equation}
{\chi^m_n}(r,s)=\frac{J_{n}(-irs)f^m_n(s)}{J_{n}(-is)}.
\end{equation}
Taking the inverse Laplace transform we obtain that
$$
\varphi^m_n(r,t)=\frac{1}{ 2\pi
i}\int_{\sigma-i\infty}^{\sigma+i\infty} \frac
{J_{n}(-irs)}{J_{n}(-is)}f^m_n(s)e^{st}\,ds,\quad n=1,2,\ldots.
$$
The right hand side of the above expression can be rewritten via
residues as in the 3D-case.

\subsubsection{Solution of IBVP for a wave equation in the interior
of a ball (by singular integral Volterra equations).}
Define
\begin{equation*}
v^m_n(x,t)=Y^m_n(\phi)\int_{t-r}^{t+r}\omega(\tau)\Psi_n\left(\frac{\tau-t}{r}\right)\,d\tau,\quad
n=1,2,\ldots,
\end{equation*}
where $\Psi_n$ are defined in (\ref{form_Psi}), and $\omega$ is
arbitrary smooth function, and check that $v_n$ satisfies wave
equation. Making the change of variables $\alpha=\frac{\tau-t}{r}$
in the integral we see that
$v^m_n(x,t)=Y^m_n(\phi)\int\limits_{-1}\limits^{1}\omega(r\alpha+t)\Psi_n(\alpha)r\,d\alpha$.
Then we evaluate:
\begin{eqnarray*}
{v^m_n}_{tt}-\Delta v^m_n={v^m_n}_{tt}-{v^m_n}_{rr}-\frac{1}{r}{v^m_n}_r-\frac{1}{r^2}{v^m_n}_{\phi\phi}=\\
=Y^m_n(\phi)\int\limits_{-1}\limits^{1}(\omega''(r\alpha+t)-\alpha^2\omega''(r\alpha+t)-\frac{3\alpha}{r}\omega'(r\alpha+t)+\\
+\frac{n^2-1}{r^2}\omega(r\alpha+t) )\Psi_n(\alpha)r\,d\alpha =\\
=Y^m_n(\phi)\int\limits_{-1}^1 \Big( \frac{1-\alpha^2}{r^2}
\frac{\partial^2}{\partial\alpha^2}\omega(r\alpha+t) -
\frac{3\alpha}{r^2} \frac{\partial}{\partial\alpha}\omega(r\alpha+t)+\\
+\frac{n^2-1}{r^2}\omega(r\alpha+t)
\Big)\Psi_n(\alpha)r\,d\alpha=\\
=\frac{Y^m_n(\phi)}{r}\int\limits_{-1}\limits^{1} \omega(r\alpha+t)\Big[ \Big((1-\alpha^2)\Psi_n(\alpha) \Big)''+3(\alpha\Psi_n(\alpha))'+(n^2-1)\Psi_n(\alpha) \Big]\,d\alpha =\\
\frac{Y^m_n(\phi)}{r}\int\limits_{-1}\limits^{1}
\omega(r\alpha+t)\Big[ (1-\alpha^2)\Psi_n''(\alpha)
-\alpha\Psi_n'(\alpha)+n^2\Psi_n(\alpha) \Big]\,d\alpha=0.
\end{eqnarray*}
The last equality is valid because of (\ref{Psi_2}). If we shift
the second argument of $v^m_n$ by $1$ we obtain that
\begin{equation}
\label{V_sol_int}
v(x,t)=\sum_{n,m}Y^m_n(\phi)\int_{t-1-r}^{t-1+r}\omega^m_n(\tau)\Psi_n\left(\frac{\tau+1-t}{r}\right)\,d\tau.
\end{equation}
Therefore (\ref{wave_eq_0}) is valid if and only if the
relation $v|_{S^1\times[0,T]}=F(t)$ holds, which is equivalent to the
equalities
$$
\int\limits_{t-2}^{t}
\omega^m_n(\tau)\Psi_n(\tau+1-t)\,d\tau=F^m_n(t),\quad n=1,2,\ldots.
$$
On differentiating the last equalities with respect to $t$ we obtain
the following equations
\begin{equation}\label{Volterra_eqn1_2d_int}
\int_{t-2}^{t}\omega^m_n(\tau)\frac{nT_n(t+1-\tau)}{\sqrt{(t+1-\tau)^2-1}}\,d\tau=(F^m_n)'(t),\quad
n=1,2,\ldots,
\end{equation}
 which are singular integral equations. Finally we come to the
\begin{proposition}
The solution of (\ref{wave_eq_0}) in 2d case can be constructed in
the form (\ref{V_sol_int}), where $\omega^m_n$ can be found as a
solutions to the singular integral equation
(\ref{Volterra_eqn1_2d_int}), with $F^m_n$ being Fourier coefficient
in the expansion (\ref{BC_expan_2}).
\end{proposition}

 \subsection{Recursive procedure for solving IP}  Consider the  Cauchy problem:
 \begin{equation}
 \label{wave_eq2}
 \left\{\begin{array}l u_{tt}-\Delta u=0,\quad x\in
 \mathbb{R}^2,\,\,t>0,\\
 u\bigl|_{t=0}=a(x),\,\,u_t\bigl|_{t=0}=b(x),
 \end{array}
 \right.
 \end{equation}
 where $\supp a,\,\supp b\subset B_1$. Let us denote $D(t,x):=\{y\in \mathbb{R}^2\,|\, |x-y|\leqslant t\}$. Due to Kirchhoff's formula, the solution
 to (\ref{wave_eq2}) admits the representation:
 $$
 u(t,x)=\frac{1}{2\pi}\int\limits_{D(t,x)} \frac{b(\xi)}{\sqrt{t^2-|\xi-x|^2}}d\xi+
 \frac{d}{dt}\left( \frac{1}{2\pi}\int\limits_{D(t,x)}
 \frac{a(\xi)}{\sqrt{t^2-|\xi-x|^2}}d\xi\right),
 $$
 and
 \begin{align*}
 u_t(t,x)=\frac{d}{dt}\left(\frac{1}{2\pi}\int\limits_{D(t,x)} \frac{b(\xi)}{\sqrt{t^2-|\xi-x|^2}}d\xi\right)\\+
 \frac{d^2}{dt^2}\left( \frac{1}{2\pi}\int\limits_{D(t,x)} \frac{a(\xi)}{\sqrt{t^2-|\xi-x|^2}}d\xi\right).
 \end{align*}
 If we take $T>2$ and change the variable $t$ to $T-t$ we come to
 $$
 a(x)=\frac{1}{2\pi}\int\limits_{D(T,x)} \frac{-u_t(T,\xi)}{\sqrt{t^2-|\xi-x|^2}}d\xi\Bigg|_{t=T}+
 \frac{d}{dt}\left( \frac{1}{2\pi}\int\limits_{D(T,x)} \frac{u(T,\xi)}{\sqrt{t^2-|\xi-x|^2}}d\xi\right)\Bigg|_{t=T}.
 $$
 Denote $K(T,x)=D(T,x)\backslash B_1$. We provide the algorithm in the case when $b=0$: in this situation
 \begin{align}
 a(x)=\frac{1}{2\pi}\int\limits_{K(T,x)} \frac{-u_t(T,\xi)}{\sqrt{t^2-|\xi-x|^2}}d\xi\Bigg|_{t=T}\label{A_form}\\+
 \frac{d}{dt}\left( \frac{1}{2\pi}\int\limits_{K(T,x)} \frac{u(T,\xi)}{\sqrt{t^2-|\xi-x|^2}}d\xi\right)\Bigg|_{t=T}\notag\\
 +\frac{1}{(2\pi)^2}\int\limits_{B_1}\int\limits_{B_1} \frac{a(\xi')}{\sqrt{T^2-|\xi-x|^2}}
 \frac{d^2}{dt^2}\frac{1}{\sqrt{t^2-|\xi-\xi'|^2}}\Bigg|_{t=T} d\xi d\xi'\notag\\
 +\frac{1}{(2\pi)^2}\int\limits_{B_1}\int\limits_{B_1} \frac{d}{dt}\frac{a(\xi')}{\sqrt{T^2-|\xi-x|^2}}\Bigg|_{t=T}
 \frac{d}{dt}\frac{1}{\sqrt{t^2-|\xi-\xi'|^2}}\Bigg|_{t=T} d\xi d\xi'\notag\\
 =I_1+I_2+I_3+I_4.\notag
 \end{align}
 Integrals $I_1$ and $I_2$ can be calculated after solving the wave equation in the exterior of a ball (\ref{wave_out}),
 using the given boundary observation $u\bigl|_{S_1\times [0,T]}=F$.
 We will show that $I_3$ and $I_4$ are decaying when $T\to\infty$.
 Denote $\beta(t,c)=\frac{1}{\sqrt{t^2-c^2}}$, then
 $$
 \beta_t(t,c)=\frac{-t}{(t^2-c^2)^{3/2}},\quad \beta_{tt}(t,c)=\frac{-1}{(t^2-c^2)^{3/2}}+\frac{3t^2}{(t^2-c^2)^{5/2}}.
 $$
We can rewrite $I_3$, $I_4$ as
\begin{align}
I_3=\frac{1}{(2\pi)^2}\int\limits_{B_1}\int\limits_{B_1} \beta(T,|\xi-x|) \beta_{tt}(T,|\xi-\xi'|) a(\xi') d\xi d\xi',\label{I3}\\
I_4=\frac{1}{(2\pi)^2}\int\limits_{B_1}\int\limits_{B_1}
\beta_t(T,|\xi-x|) \beta_t(T,|\xi-\xi'|) a(\xi') d\xi
d\xi'.\label{I4}
\end{align}
When $|\xi|<1$, $|\xi'|<1$ and $|x|<1$, then
\begin{equation}
\label{Beta_ner}
|\beta|<\frac{1}{T^2-2^2},\quad
|\beta_t|<\frac{T}{(T^2-2^2)^{3/2}},\quad
|\beta_{tt}|<\frac{4T}{(T^2-2^2)^{3}}.
\end{equation}
Therefore $I_3$ and $I_4$ can be made arbitrarily small taking
time of observation $T\gg 2$. If we restrict ourselves to smaller
$T>2$, we can provide the iterative procedure of recovering
initial data $a$. On introducing the operator $K:C(B_1)\mapsto
C(B_1)$, acting by the rule $Ka:=I_3a+I_4a,$ the equality
(\ref{A_form}) can be written in a form
\begin{equation}
\label{Iter_eqn} a=I_1+I_2+I_3+I_4=CF+Ka=a_1+Ka.
\end{equation}
The latter equation should be understood in the following way:
initial displacement $a$ is equal to sum of $a_1$, the term which
can be calculated from observation $F$, and additional term $Ka$.
Note that the norm of operator $K$ can be made small by the choice
of $T$ (see (\ref{I3}), (\ref{I4}) and (\ref{Beta_ner})). On a
next step one takes $a_1$ as initial data and solve the forward
problem (\ref{wave_eq}) in the case $n=2$ with $a=a_1,$ $b=0$ to
find the new observation $u\bigl|_{S_1\times [0,T]}=F_1$. Then one
need to repeat the procedure described above with the observation
$F-F_1$ and find $a_2=C(F-F_1)$. As before
$$
a-a_1=C(F-F_1)+K(a-a_1)=a_2+K(a-a_1),
$$
and therefore
$$
a-(a_1+a_2)=K(a-a_1)=K^2a.
$$
We can continue this process, solving the forward problem for
(\ref{wave_eq}) with $a=a_k$, $b=0$ to find
$F_k=u\bigl|_{S_1\times [0,T]}$. The repetition  of this procedure
$n$ times gives the following approximation formula:
$$
a=\sum_{i=1}^n a_i+K^n a,
$$
where the second term in the right hand side is small by the
choice of $T$ such that $\|K\|<1$. Thus  initial data $a$ are
approximated by the sum of $a_i$, where each $a_i$ is a solution
to IP constructed by our algorithm.

\section{Remarks.}
It is interesting to compare our results with ones of Rakesh and
D. Finch. In \cite{FR} they prove several striking inverse
formulas, one of them has the following form (in 3d case and if
displacement  $a\equiv 0$ in (\ref{wave_eq}))
\begin{equation}\label{FRFormula}
b(x)=-\frac{1}{2\pi}\int\limits_{S_1}\frac{\partial^2_t(tF(t,x_0))|_{t=|x-x_0|}}{|x-x_0|}\,dS_1(x_0).
\end{equation}
We note, that for reconstruction of initial speed $b$ at point $x$
one need to know boundary measurement $F(\cdot,t)$ on times up to
$|x|+1$ and differentiate function $F$ twice with respect to $t$.
Our reconstruction procedure makes almost the same, but in terms
of the Fourier coefficients (it is also necessary to differentiate
functions $f_n$ twice with respect to $t$ and to find initial
conditions $a$ and $b$ at point $x$ we need to know boundary
measurements $f_n(\cdot,t$) on the same time interval up to
$T=|x|+1$.  Therefore our procedure can be considered as a
discrete analog of Rakesh and Finch's formula. The peculiarity of
our procedure is that using the observation in the time interval
$(0,|x|+1)$, one can restore both initial conditions $a$ and $b$
from (\ref{wave_eq}) at point $x$.

On the other hand if one of the Cauchy data $a$ and $b$ are equal to
zero then a simple modification of our procedure enables one to
restore the remaining function using the observation on half time
interval. For example, let $a=0$ in (\ref{wave_eq}), i.e. $u$ is
the solution to the following  Cauchy problem:
\begin{equation*}
\left\{\begin{array}l u_{tt}-\Delta u=0,\quad x\in
\mathbb{R}^3,\,(\text{or $x\in \mathbb{R}^2$}),\,\,t>0,\\
u\bigl|_{t=0}=0,\,\,u_t\bigl|_{t=0}=b(x),
\end{array}
\right.
\end{equation*}
where $b\in C^2_0(B_1)$ and observation
$F:=u\bigl|_{S_1\times(0,T)}$ is given with $T=1$. Note that in
this case reconstruction formula (\ref{FRFormula}) is not
applicable.

On introducing the odd reflection of $u(x,t)$ with respect to
$t=0$ by the rule:
\begin{equation*}
U(x,t)=\left\{\begin{array}l u(x,t),\quad t>0,\\
-u(x,t),\quad t<0,
\end{array}
\right.
\end{equation*}
we see that $U(x,t)$ satisfies
\begin{equation*}
U_{tt}-\Delta U=0,\quad x\in \mathbb{R}^3,\,\,t\in \mathbb{R},
\end{equation*}
and new observation
\begin{equation*}
\widetilde F:= U\bigl|_{S\times \mathbb{R}}=\left\{\begin{array}l F(x,t),\quad t>0,\\
-F(x,t),\quad t<0
\end{array}
\right.
\end{equation*}
is known on $S\times (-1,1)$. Now we can introduce the function
\begin{equation*}
V(x,t):=U(x,t-1), \quad x\in \mathbb{R}^3,
\end{equation*}
then for $t>0$ this function satisfies the wave equation
\begin{equation*}
V_{tt}-\Delta V=0,\quad x\in \mathbb{R}^3,\,\, 0<t<2,
\end{equation*}
and we know $V$ on $S\times (0,2)$, where $V\bigl|_{S\times
(0,2)}=\widetilde F(\cdot,\cdot -1)\bigl|_{S\times (0,2)}$. Then
we can apply the procedure described in the second section if we are
in 3d case or in the third section if we are in 2d case to recover
$V$ and $V_t$ at ${t=0}$. It remains to use Kirchhoff formula to
obtain initial condition $b(x)=V(x,1)$.

We finalize the above arguments in the following remark:
\begin{remark}
If in the IBVP (\ref{wave_eq}) in $R^n,$ $n=\{2,3\}$ one of the Cauchy
data is zero, then we can recover the remaining one using the
observation (\ref{observation}) on $S\times (0,1)$.
\end{remark}

\noindent{\bf Acknowledgments}

The research of Mikhail Belishev  was supported in part by RFBR 17-01-00529a.
Alexandr Mikhaylov and Victor Mikhaylov were supported by RFBR 18-01-00269.
This work is supported by the Program of the Presidium of the
Russian Academy of Sciences 01 Fundamental Mathematics and
its Applications under the grant PRAS-18-01.
We thank the Volkswagen Foundation (VolkswagenStiftung)  program ``Modeling, Analysis, and Approximation Theory
toward application in tomography and inverse problems'' for kind support and
stimulating our collaboration.

\end{document}